\documentclass[twocolumn]{article}       
\RequirePackage{fix-cm}
\usepackage{graphicx}
\usepackage[english]{babel}
\usepackage[utf8]{inputenc}
\usepackage{todonotes}
\usepackage{amsmath,amssymb}
\usepackage{natbib}
\usepackage{float}
\usepackage{ulem}
\usepackage{algorithm,algorithmicx}
\usepackage{algpseudocode}
\usepackage{url}

\def\RR{\textsf{R}\/}

\let\pkg=\strong

\newcommand{\x}{\mathbf{x}}

\newcommand{\bolda}{\mathbf{a}}
\newcommand{\boldalpha}{\boldsymbol{\alpha}}
\newcommand{\estboldalpha}{\hat{\boldsymbol{\alpha}}}
\newcommand{\estb}{\hat{\mathbf{b}}}
\newcommand{\estalpha}{\hat{\boldsymbol{\alpha}}}

\newcommand{\X}{\mathbf{X}}
\newcommand{\estSigma}{\widehat{\Sigma}}
\newcommand{\estGamma}{\widehat{\Gamma}}
\newcommand{\slice}[1]{\mathcal{S}_{#1}}

\newcommand{\meanX}{\overline{X}}
\newcommand{\R}[1]{\mathbb{R}^{#1}}
\newcommand{\proj}{\mathcal{P}}
\newcommand{\boldproj}{\mathbf{P}}
\newcommand{\estPi}{\widehat{\Pi}}

\begin{document}
\sloppy

\title{Interpretable sparse SIR for functional data}

\author{Victor Picheny \and Rémi Servien \and Nathalie Villa-Vialaneix}


\date{Received: date / Accepted: date}

\maketitle

\begin{abstract}
We propose a semiparametric framework based on Sliced Inverse Regression (SIR) to address 
the issue of variable selection in functional regression. SIR is an effective method 
for dimension reduction which computes a linear 
projection of the predictors in a low-dimensional space, without loss of 
information on the 
regression. In order to deal with the high dimensionality of the predictors, 
we consider penalized versions of SIR: ridge and sparse.
We extend the approaches of variable selection developed 
for multidimensional SIR to select intervals that form a partition of the 
definition domain of the functional predictors. Selecting entire intervals 
rather than separated evaluation points improves the interpretability of 
the estimated coefficients in the 
functional framework.
A fully automated iterative procedure is 
proposed to find the critical (interpretable) intervals. The approach is 
proved efficient on simulated and real data. The method is implemented in the 
\RR{} package \pkg{SISIR} available on CRAN at 
\url{https://cran.r-project.org/package=SISIR}.
\end{abstract}

\section{Introduction}

This article focuses on the functional regression problem, in which a real 
random variable $Y$ is predicted from a functional predictor $X(t)$ that takes 
values in a functional space ({\it e.g.}, $L^2([0,1])$, the space of squared 
integrable functions over $[0,1]$), based on a set of observed pairs $(X,Y)$, 
$(\x_i,y_i)_{i=1,\ldots,n}$. The main challenge with functional regression 
lies in its high dimension: the underlying dimension of a functional space is 
infinite, and even if the digitized version of the curves is considered, the 
number of evaluation points is typically much larger than the number of 
observations.

Recently, an increasing number of works have focused on variable selection in 
this functional regression framework, in particular in the linear
setting. The problem is to select parts of the 
definition domain of $X$ that are relevant to predict $Y$. Considering digitized 
versions of the functional predictor $X$, approaches based on Lasso have been 
proposed to select a few isolated points of $X$ 
\citep{ferraty_etal_B2010,aneiros_vieu_SPL2014,mckeague_sen_AS2010,
kneip_etal_AS2016}. Alternatively, other authors proposed to perform variable 
selection on predefined functional bases. For instance, 
\citet{matsui_konishi_CSDA2011} used $L^1$ regularization on Gaussian basis 
functions and \citet{zhao_etal_JCGS2012,chen_etal_SIAMJSC2015} on wavelets.

However, in many practical situations, the relevant information may not correspond 
to isolated evaluation points of $X$ neither to some of the components of its expansion on 
a functional basis, but to its value on some continuous intervals, $X(\left[t_a, t_b\right])$. 
In that case, variable selection amounts to identify those intervals. As advocated by 
\citet{james_etal_AS2009}, a desirable feature of variable selection provided by 
such an approach is to enhance the interpretability of the relation between $X$ 
and $Y$. Indeed, it reduces the definition domain of the 
predictors to a few influential intervals, or it focuses on some particular 
aspects 
of the curves in order to obtain expected values for $Y$. 
Tackling this issue can be seen as selecting groups of contiguous variables 
({\it i.e.}, intervals) instead of selecting isolated variables. 
\citet{fraiman_etal_JMA2016}, in the linear setting, and 
\citet{fauvel_etal_IEEEJSTAEORS2015,ferraty_hall_JCGS2015}, in a nonparametric 
framework, propose several alternatives to do so. However, no specific 
contiguity constraint is put on groups of variables.

In the present work, we propose a semi-parametric model that selects intervals 
in the definition domain of $X$ with an automatic approach. The method is based 
on Sliced Inverse Regression \citep[SIR,][]{li_JASA1991}: the main idea of SIR 
is to define a low dimensional data-driven subspace on which the functional 
predictors can be projected. This subspace, called Effective Dimension 
Reduction (EDR) space is defined so as to optimize the prediction ability of 
the projection. As a particular case, the method includes the linear 
regression. Our choice for SIR is motivated by the fact that the method is based 
on a semi-parametric model that is more flexible than linear models. The method 
has been extended to the functional framework in previous works
\citep{ferre_yao_S2003,ferre_villa_SJS2006} and sparse ({\it i.e.}, $\ell_1$
penalized) versions of the approach have also already been proposed in 
\citet{li_nachtsheim_T2008} and \citet{li_yin_B2008} for the multivariate 
framework. Building on these previous proposals, we show that a tailored 
group-Lasso-like penalty allows us to select groups of variables corresponding 
to intervals in the definition domain of the functional predictors. 

Our second contribution is a fast and automatic procedure for building intervals 
in the definition domain of the predictors without using any prior knowledge.
As far as we know, the only works that propose a method to both define 
and select relevant intervals in the domain of the predictors are the work of
\citet{park_etal_p2016} and \citet{grollemund_etal_p2018}, both in the linear 
framework.
Our approach is based 
on an iterative procedure that uses the full regularization path of the Lasso.

The paper is organized as follows: Section~\ref{reg-sir} presents the SIR 
approach in a multidimensional framework and its adaptations to the 
high-dimensional and functional frameworks, which are based on regularization 
and/or sparsity constraints. Section~\ref{il-sir} describes our proposal when
the domain of the predictors are partitioned using a fixed set of intervals. 
Then, Section~\ref{interv} describes an automatic procedure to find these 
intervals and Section~\ref{param-selection} provides practical methods to tune 
the different parameters in a high dimensional framework. Finally, 
Section~\ref{simulations} evaluates our approach on simulated and real-world 
datasets.

\section{A review on SIR and regularized versions}\label{reg-sir}

In this section, we review the standard SIR for multivariate data and its 
extensions to the high-dimensional setting. Here, $(X,Y)$ denotes a random pair 
of variables such that $X$ takes values in $\R{p}$ and $Y$ is real. We assume 
given $n$ i.i.d. realizations of $(X,Y)$, $(\x_i,y_i)_{i=1,\ldots,n}$.

\subsection{The standard multidimensional case}

When $p$ is large, classical modeling approaches suffer from the curse of 
dimensionality. 
This problem might occur even if $p$ is smaller than $n$. A standard way to 
overcome this issue is to rely on dimension reduction techniques. This kind of 
approaches is based on the assumption that there exists an Effective 
Dimension Reduction (EDR) space $\mathcal{S}_{Y|X}$ which is the smallest 
subspace such that the projection of $X$ on $\mathcal{S}_{Y|X}$ retains all the 
information on $Y$ contained in the predictor $X$. More precisely, 
$\mathcal{S}_{Y|X}$ is assumed of the form 
$\mbox{Span}\{\bolda_1,\ldots,\bolda_d\}$, with $d \ll p$, such that 
\begin{equation}
	\label{eq::sir-model}
	Y = F(\bolda_1^\top X, \ldots, \bolda_d^\top X, \epsilon),
\end{equation}
in which $F:\mathbb{R}^{p+1}\rightarrow\mathbb{R}$ is an unknown function and 
$\epsilon$ is an error term independent of $X$. To estimate this subspace, SIR 
is one of the most classical approaches when $p<n$: under an appropriate and 
general 
enough condition, \citet{li_JASA1991} shows that $\bolda_1,\ldots,\bolda_d$ can 
be estimated as the first $d$ $\Sigma$-orthonormal eigenvectors of the 
generalized eigenvalue problem: $\Gamma \bolda = \lambda \Sigma \bolda$, in 
which $\Sigma$ is the covariance matrix of $X$ and $\Gamma$ is the covariance 
matrix of $\mathbb{E}(X|Y)$. 

In practice, $\Sigma$ is replaced by the empirical covariance, $\estSigma = 
\frac{1}{n} \sum_{i=1}^n \left(\x_i - \meanX\right) \left(\x_i- 
\meanX\right)^\top$, and $\Gamma$ is estimated by ``slicing'' the observations 
$(y_i)_i$ as follows. The range of $Y$ is partitioned into $H$ consecutive and 
non-overlapping slices, denoted hereafter $\mathcal{S}_1$, \ldots, 
$\mathcal{S}_H$. An estimate of $\mathbb{E}(X|Y)$ is thus simply obtained by 
$\left(\meanX_1, \ldots,\meanX_H\right)$ in which $\meanX_h$ is the average of 
the observations $\x_i$ such that $y_i$ is in $\slice{h}$ and $\meanX_h$ is 
associated with the empirical frequency $\hat{p}_h = \frac{n_h}{n}$ with $n_h$ 
the 
number of observations in $\slice{h}$. $\estGamma$ is thus defined as 
$\sum_{h=1}^H \hat{p}_h \meanX_h \meanX_h^\top$.

SIR has different equivalent formulations that can be useful to introduce 
regularization and sparsity.
\citet{cook_AS2004} 
shows that the SIR estimate can be obtained by minimizing 
over $A\in \mathbb{R}^{p\times d}$ and $C = \left(C_1,\,...,\,C_H\right)$, 
with $C_h \in \R{d}$ (for $h=1,\ldots,H$),
\begin{equation}
	\label{eq::minpb-sir}
	\mathcal{E}_1(A,C) = \sum_{h=1}^H \hat{p}_h \|(\meanX_h - \meanX) - \estSigma 
AC_h\|^2_{\estSigma^{-1}},
\end{equation}
in which $\|.\|_{\estSigma^{-1}}^2$ is the norm $\forall\,u\in\R{p}$, 
$\|u\|_{\estSigma^{-1}}^2 = u^\top \estSigma^{-1} u$ and the searched vectors 
$\bolda_j$ are the columns of $A$.

An alternative formulation is described in \citet{chen_li_SS1998}, where SIR
is written as the following optimization problem:
\begin{equation}
	\label{eq::opt-sir}
	\max_{\bolda_j, \phi} \mbox{Cor}(\phi(Y), \bolda_j^\top X),
\end{equation}
where $\phi$ is any function $\R{} \rightarrow \R{}$ and $(\bolda_j)_j$ are 
$\Sigma$-orthonormal. So, SIR can be interpreted as a canonical 
correlation problem. The authors also prove that the solution of $\phi$ 
optimizing Equation~(\ref{eq::opt-sir}) for a given $\bolda_j$ is $\phi(y) = 
\bolda_j^\top \mathbb{E}(X|Y=y)$, and that $\bolda_j$ is also obtained as the 
solution of the mean square error optimization $\min_{\bolda_j} 
\mathbb{E}\left(\phi(Y) - \bolda_j^\top X\right)^2$.

However, as explained in \citet{li_yin_B2008} and 
\citet{coudret_etal_JSFdS2014} among others, in a high dimensional setting ($n < 
p$), $\estSigma$ is singular and the SIR problem is thus ill-posed. The same 
problem occurs in the functional setting \citep{dauxois_etal_CRAS2001}. 
Solutions to overcome this difficulty include variable selection 
\citep{coudret_etal_JSFdS2014}, ridge regularization or sparsity constraints.

\subsection{Regularization in the high-dimensional setting}\label{ridge-sir}

In the high-dimensional setting, directly applying a ridge penalty, $ \mu_2 
\sum_{h=1}^H \hat{p}_h \|AC_h\|_{\mathbb{I}_p}^2$ (for a given $\mu_2 > 0$), to 
$\mathcal{E}_1$ would require the computation of $\estSigma^{-1}$ (see 
Equation~(\ref{eq::minpb-sir})) that does not exist when $n < p$. 
However, \cite{bernardmichel_etal_B2008} show that this problem can 
be rewritten as the minimization of
\begin{eqnarray}
	\label{eq::ridge-sir-2}
	&&\sum_{h=1}^H \hat{p}_h C_h^\top A^\top (\estSigma + \mu_2 \mathbb{I}_p) A 
C_h -\nonumber\\
	&&\qquad \qquad \qquad 2 \sum_{h=1}^H \hat{p}_h \left(\meanX_h - \meanX 
\right) A C_h,
\end{eqnarray}
which is well defined even for the high-dimensional setting. Minimizing this 
quantity with respect to $A$ leads to define the columns of $A$ (and hence the 
searched vectors $\bolda_j$) as the first $d$ eigenvectors of $\left(\estSigma 
+ \mu_2 \mathbb{I}_p \right)^{-1} \estGamma$.

\subsection{Sparse SIR}

Sparse estimates of $\bolda_j$ usually increase the interpretability of the 
model (here, of the EDR space) by focusing on the most important predictors 
only. Also, \cite{lin_etal_AS2018} prove the relevance of sparsity for SIR in 
high dimensional setting by proposing a consistent screening pre-processing of 
the variables before the SIR estimation. A different and very commmon approach 
is to handle sparsity directly by a sparse penalty (in the line of the 
well-known Lasso). However, contrary to ridge regression, adding directly a 
sparse penalty to Equation~(\ref{eq::minpb-sir}) does not allow a reformulation 
valid for the case $n < p$. To the best of our knowledge, only two alternatives 
have already been published to use such methods, one based on the regression 
formulation (\ref{eq::minpb-sir}) and the other on the correlation formulation 
(\ref{eq::opt-sir}) of SIR. 

\cite{li_yin_B2008} derive a sparse ridge estimator from the work of 
\cite{ni_etal_B2005}. Given $(\hat{A},\hat{C})$, solution of the ridge SIR, a 
shrinkage index vector $\boldalpha = (\alpha_1, \ldots, \alpha_p)^\top \in 
\R{p}$ is obtained by minimizing a least square error with $\ell_1$ penalty: 
\begin{eqnarray}
	\label{eq::liyin_sparse}
\mathcal{E}_{s,1}(\boldalpha)& =& \sum_{h=1}^H \hat{p}_h \left\| \left( \meanX_h 
- \meanX \right) - \estSigma \mbox{Diag}(\boldalpha) \hat{A}\hat{C}_h \right\|_{\mathbb{I}_p}^2 +\nonumber\\
&& \qquad \mu_1 \|\boldalpha\|_{\ell_1}, 
\end{eqnarray}
 for a given $\mu_1 \in \R{+*}$ where $\|\boldalpha\|_{\ell_1}=\sum_{j=1}^p 
|\alpha_j|$. Once the coefficients $\boldalpha$ have been estimated, the EDR 
space is the space spanned by the columns of $\mbox{Diag}(\estalpha) \hat{A}$, 
where $\estalpha$ is the solution of the minimization of 
$\mathcal{E}_{s,1}(\boldalpha)$.

An alternative is described in 
\cite{li_nachtsheim_T2008} using the
correlation formulation of the SIR. After the standard SIR estimates $\hat{\bolda}_1,\ldots,\hat{\bolda}_d$ have 
been computed, they solve $d$ independent minimization 
problems with sparsity constraints introduced as an $\ell_1$ penalty: 
$\forall\,j=1, \ldots, d$,
\begin{equation}
	\label{eq::sir-mlr}
	\mathcal{E}_{s,2}(\bolda^s_j) = \sum_{i=1}^n 
\left[\mathcal{P}_{\hat{\bolda}_j}(X|y_i) - (\bolda_j^s)^\top \x_i\right]^2 + 
\mu_{1,j} \|\bolda_j^s\|_{\ell_1},
\end{equation}
in which $\mathcal{P}_{\hat{\bolda}_j}(X|y_i) = 
\widehat{\mathbb{E}}(X|Y=y_i)^\top 
\hat{\bolda}^{j}$,  with $\widehat{\mathbb{E}}(X|Y=y_i) = \meanX_h$ for $h$ 
such 
that $y_i \in \slice{h}$ in the case of a sliced estimate of 
$\widehat{\mathbb{E}}(X|Y)$ and $\mu_{1,j}>0$ is a parameter controlling the 
sparsity of the solution.

Note that both proposals have problems in the high-dimensional setting:
\begin{itemize}
	\item In their proposal, \cite{li_yin_B2008} avoid the issue of the 
singularity of $\estSigma$ by working in the original scale of the predictors 
for both the ridge and the sparse approach (hence the use of the 
$\|.\|_{\mathbb{I}_p}$-norm in Equation~(\ref{eq::liyin_sparse}) instead of the 
standard $\|.\|_{\estSigma^{-1}}$-norm of Equation~(\ref{eq::minpb-sir})).
However, for the ridge problem, this choice has been proved to produce a 
degenerate problem by \cite{bernardmichel_etal_B2008}.
	\item \cite{li_nachtsheim_T2008} base their sparse version of the SIR on the 
standard estimates of the SIR problem that cannot be computed in the 
high-dimensional setting. 
\end{itemize}

Moreover, the other differences between these two approaches can be summarized 
in two points:
\begin{itemize}
	\item using the approach of \cite{li_yin_B2008} based on shrinkage
coefficients, the indices $\alpha_j$ where $\alpha_j>0$ are the same on all the 
$d$ components of the EDR. This makes sense because the vectors $\bolda_j$ 
themselves are not relevant: only the space spanned by them is and so there is 
no interest to select different variables $j$ for the $d$ estimated 
directions. 
Moreover, this allows to formulate the optimization in a single problem. 
However, this problem relies on a least square minimization with dependent 
variables in a high dimensional space $\R{p}$;
	\item on the contrary, the approach of \cite{chen_li_SS1998} relies on
a least square problem based on projections and is thus obtained from $d$ 
independent optimization problems. The dimension of the dependent variable is 
reduced (1 instead of $p$) but the different vectors 
which span the EDR space are estimated independently and not simultaneously.
\end{itemize}

In our proposal, we combine both advantages of \cite{li_yin_B2008} and \cite{li_nachtsheim_T2008} using a single optimization problem 
based on the correlation formulation of SIR. In this problem, the dimension of 
the dependent variable is reduced ($d$ instead of $p$) when compared to the 
approach of 
\cite{li_yin_B2008} and it is thus computationally more efficient.
Identical sparsity 
constraints are imposed on all $d$ dimensions using a shrinkage approach, but 
instead 
of selecting the nonzero variables independently, we adapt the sparsity 
constraint to the functional setting to avoid selecting isolated measurement 
points. The next section describes this approach.

\section{Sparse and Interpretable SIR (SISIR)}\label{il-sir}

In this section, a functional regression framework is assumed. $X$ is thus a 
functional random variable, taking value in a (infinite dimensional) Hilbert 
space. $(x_i,y_i)_{i=1,\ldots,n}$ are $n$ i.i.d. realizations of $(X,Y)$. 
However, $x_i$ are not perfectly known but observed on a given (deterministic) 
grid $\tau = \{t_1, \ldots, t_p\}$. We denote by $\x_i = 
(x_i(t_j))_{j=1,\ldots,p} \in \R{p}$ the $i$-th observation, by $\x^j = 
(x_i(t_j))_{i=1,\ldots,n}$ the observations at $t_j$ and by $\X$ the $n 
\times p$ matrix $(\x_1, \ldots, \x_n)^\top$. Unless said otherwise, the 
notations 
are derived from the ones introduced in the multidimensional setting 
(Section~\ref{reg-sir}) by using the $\x_i$ as realizations of $X$.

Some very common methods in functional data analysis, such as 
splines \citep{hastie_etal_ESL2001}, use the supposed smoothness of $X$ to 
project 
them in a reduced dimension space. 
Contrary to these methods, we do not use or need that the observed functional 
predictor is smooth. We take advantage of the functional 
aspects of the data in a different way, using the natural ordering of the 
definition domain of $X$ to impose sparsity on the EDR space. To do so, 
we assume that this definition domain is partitioned into $D$ contiguous and 
non-overlapping intervals, $\tau_1, \ldots, \tau_D$. In the present section, 
these intervals are supposed to be given {\it a priori} and we will describe 
later (in Section~\ref{interv}) a fully automated procedure to obtain them from 
the data.

The following two subsections are devoted to the description of the two steps 
(ridge and sparse) of the method, adapted from
\cite{bernardmichel_etal_B2008,li_yin_B2008,li_nachtsheim_T2008}.

\subsection{Ridge estimation}\label{ilsir-ridge}

The ridge step is the minimization of Equation~(\ref{eq::ridge-sir-2}), over 
$(A,C)$ to obtain $\hat{A}$ and $\hat{C}$. In practice, the solution is computed 
as follows:

\begin{enumerate}
	\item The estimator of $A\in\mathbb{R}^{p \times d}$ is the solution of the 
ridge penalized SIR and is composed of the first $d$ $\left(\estSigma + 
\mu_2 \mathbb{I}_p\right)$-orthonormal eigenvectors of $\left(\estSigma + \mu_2 
\mathbb{I}_p\right)^{-1} \estGamma$ associated with the $d$ largest eigenvalues. 
In practice, the same procedure as the one described in 
\cite{ferre_yao_S2003,ferre_villa_SJS2006} is used: first, orthonormal 
eigenvectors (denoted hereafter $(\estb_j)_{j=1,\ldots,d}$) of the matrix 
$\left(\estSigma + \mu_2 \mathbb{I}_p\right)^{-1/2} \estGamma \left(\estSigma + 
\mu_2 \mathbb{I}_p\right)^{-1/2}$ are computed. Then, $\hat{A}$ is the matrix 
whose columns are equal to $\left(\estSigma + \mu_2 \mathbb{I}_p\right)^{-1/2} 
\estb_j$ for $j=1,\ldots,d$. It is easy to prove that these columns are 
$\left(\estSigma + \mu_2 \mathbb{I}_p\right)$-orthonormal eigenvectors of 
$\left(\estSigma + \mu_2 \mathbb{I}_p\right)^{-1} \estGamma$.

	\item For a given $A$, the optimal $\hat{C} = (\hat{C}_1, \ldots, \hat{C}_H) 
\in \mathbb{R}^{d,H}$ is given by the first order derivation condition over the 
minimized criterion. This is equivalent to $\left[A^\top 
\estSigma A + \mu_2 A^\top A \right] 
\hat{C}_h = A^\top \left(\meanX_h - \meanX \right)$ that gives $\hat C_h = 
\left[ A^\top  \widehat \Sigma A + \mu_2 A^\top A \right]^{-1} A^\top \left( 
\meanX_h - \meanX \right) = A^\top \left( \meanX_h - \meanX \right)$ because the 
columns of $A$ are $\left(\estSigma + \mu_2 \mathbb{I}_d\right)$-orthonormal.
\end{enumerate}

\subsection{Interval-sparse estimation}\label{sisir-sparse}

Once $\hat{A}$ and $\hat C$ have been computed, the estimated projections of 
$(\widehat{\mathbb{E}}(X|Y=y_i))_{i=1,\ldots,n}$ onto the EDR space are obtained 
by: $\mathcal{P}_{\hat A}(\widehat{\mathbb{E}}(X|Y=y_i)) = 
(\meanX_h-\meanX)^\top\hat{A}$, for $h$ such that $y_i\in \slice{h}$. This $d$ 
dimensional vector will be denoted by $(\mathcal{P}_i^1, \ldots, 
\mathcal{P}_i^d)^\top$. In addition, we will also denote by $\boldproj^j$ (for 
$j=1,\ldots,d$), $\boldproj^j = (\proj_1^j, \ldots, \proj_n^j)^\top \in \R{n}$.

$D$ shrinkage coefficients, $\boldalpha = (\alpha_1, \ldots, \alpha_D) \in 
\R{D}$, one for each interval $(\tau_k)_{k=1,\ldots,D}$, are finally estimated.
If $\Lambda(\boldalpha) = \mbox{Diag} \left(\alpha_1 
\mathbb{I}_{|\tau_1|},\ldots,\alpha_D \mathbb{I}_{|\tau_D|}\right) \in 
\mathbb{R}^{p\times p}$, this leads to solve the following Lasso problem
\begin{eqnarray*}
	\estalpha &=& \arg\min_{\boldalpha \in \R{D}} \sum_{j=1}^d \sum_{i=1}^n 
\|\proj_i^j - \left(\Lambda(\boldalpha)\, \hat{\bolda}_j\right)^\top \x_i 
\|^2 + \mu_{1} \|\boldalpha\|_{\ell_1}\\
	&=& \arg\min_{\boldalpha \in \R{D}} \sum_{j=1}^d \|\boldproj^j - \left(\X 
\Delta(\hat{\bolda}_j)\right) \boldalpha \|^2 + \mu_{1} \|\boldalpha\|_{\ell_1},
\end{eqnarray*}
with $\Delta(\hat{\bolda}_j)$ the $(p\times D)$-matrix such that 
$\Delta_{lk}(\hat{\bolda}_j)$, is the $l$-th entry of $\hat{\bolda}_j$, 
$\hat{a}_{jl}$, if $t_l \in \tau_k$ and 0 otherwise. 

This problem can be rewritten as
\begin{equation}
	\label{eq::sparse-optim}
	\arg\min_{\boldalpha \in \R{D}} \|\boldproj - \Delta(\X \hat{A})\,
\boldalpha \|^2 + \mu_{1} \|\boldalpha\|_{\ell_1}
\end{equation}
with $\boldproj = \left(\begin{array}{c}
	\boldproj^1\\
	\vdots\\
	\boldproj^d
\end{array}\right)$, a vector of size $dn$ and $\Delta(\X \hat{A}) = 
\left(\begin{array}{c}
	\X \Delta(\hat{\bolda}_1)\\
	\vdots\\
	\X \Delta(\hat{\bolda}_p)\\
\end{array}\right)$, a $(dn)\times D$-matrix.

$\hat{\boldalpha}$ are used to define the $\hat{\bolda}_j^s$ of the vectors 
spanning the EDR space by:
\begin{equation*}
	\forall\,l=1,\ldots,p,\ \hat{a}_{jl}^s = \hat{\alpha}_k\, \hat{a}_{jl} \ 
\mbox{for }k\mbox{ such that }t_l \in \tau_k.
\end{equation*}

Once the sparse vectors $(\hat{\bolda}_j^s)_{j=1,\ldots,d}$ have been obtained, 
an Hilbert-Schmidt orthonormalization approach is used to make them 
$\estSigma$-orthonormal. 

Of note, as a single shrinkage coefficient is defined for all 
$(\hat{a}_{jl})_{t_l \in \tau_k}$, the method is close to group-Lasso 
\citep{simon_etal_JCGS2013}, in the sense that, for a given $k\in 
\{1,\ldots,D\}$, estimated $(\hat{a}^s_{jl})_{j=1,\ldots,d,\,t_l \in \tau_k}$ 
are either all zero or either all different from zero. However, the approach 
differs from group-Lasso because group-sparsity is not controlled by the 
$L_2$-norm of the group but by a single shrinkage coefficient associated to 
that group: the final optimization problem of Equation~(\ref{eq::sparse-optim}) 
is thus written as a standard Lasso problem (on $\boldalpha$).

Another alternative would have been to use fused-Lasso \citep{tibshirani_etal_JRSSB2005} 
to control the total variation norm of the estimates. However, the method does 
not explicitely select intervals and, as illustrated in 
Section~\ref{simulated}, is better designed to produce piecewise constant 
solutions than solutions that have sparsity properties on intervals of the 
definition domain.

\section{An iterative procedure to select the intervals}
\label{interv}

The previous section described our proposal to detect the subset of relevant 
intervals among a fixed, predefined set of intervals of the definition 
domain of the predictor, $(\tau_k)_{k=1,\ldots,D}$. However, choosing \textit{a 
priori} a proper set of intervals is a challenging task without expert 
knowledge, and a poor choice (too small, too large, or shifted intervals) 
may largely hinder interpretability. In the present section, we propose an 
iterative method to automatically design the intervals, without making any {\it 
a priori} choice. 

In a closely related framework, \cite{fruth_etal_RESS2015} tackle the problem 
of designing intervals by combining sensitivity indices, linear regression 
models and a method called sequential bifurcation 
\citep{bettonvil_CISSC1995} which allows them to sequentially split in two the
most promising intervals (starting from a unique interval covering the entire 
domain of $X$). Here, we propose the inverse approach: we start with small 
intervals and merge them sequentially. Our approach is based on the standard 
sparse SIR (which is used as a starting point) and iteratively performs the 
most relevant merges in a flexible way (contrary to a splitting approach, we do 
not need to arbitrary set the splitting positions).

The intervals $(\tau_k)_{k=1,\ldots,D}$ 
are first initialized to a very fine grid, taking for instance $\tau_k = 
\{t_k\}$ 
for all $k=1,\ldots,p$ (hence, at the beginning of the procedure, $D=p$). The 
sparse step described in Section~\ref{sisir-sparse} is then performed with the 
{\it a priori} intervals $(\tau_k)_{k=1,\ldots,D}$: the set of solutions of 
Equation~(\ref{eq::sparse-optim}), for varying values of the regularization 
parameter $\mu_1$, is obtained using a regularization path approach, as 
described in \cite{friedman_etal_JSS2010}. Three elements are retrieved from 
the 
path results:
\begin{itemize}
	\item $(\estboldalpha_k^*)_{k=1,\ldots,D}$ are the solutions of the sparse 
problem for the value $\mu_1^*$ of $\mu_1$ that minimizes the GCV error;
	\item $(\estboldalpha_k^+)_{k=1,\ldots,D}$ and 
$(\estboldalpha_k^-)_{k=1,\ldots,D}$ are the first solutions, among the path of 
solutions, such that at most (resp. at least) a proportion $P$ of the 
coefficients are non zero coefficients (resp. are zero coefficients), for a 
given chosen $P$, which should be small (0.05 for instance).
\end{itemize}

\begin{figure}[ht]
	\centering
	\includegraphics[width=\linewidth]{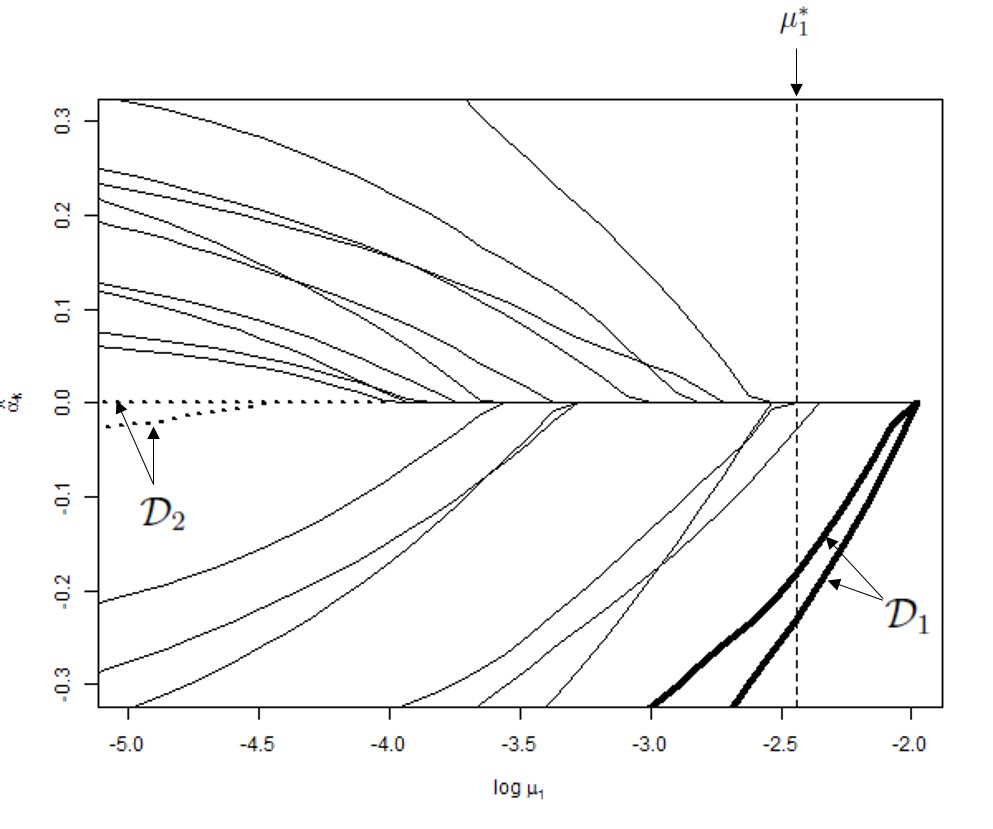}
	\caption{Example of regularization path with $D=20 $: 
$(\estboldalpha_k)_{k=1,\ldots,D}$ are plotted according to different values of 
the tuning parameter $\mu_1$. The vertical dotted line represents the optimal 
value $\mu_1^*$ that provides the solutions 
$(\estboldalpha_k^*)_{k=1,\ldots,D}$ of the sparse problem.  
$(\estboldalpha_k)_{k \in \mathcal{D}_1}$ and $(\estboldalpha_k)_{k \in 
\mathcal{D}_2}$ are respectively represented in bold and in pointed lines for 
$P=0.1$.}
	\label{fig::path}
\end{figure}

\begin{figure*}[ht]
	\centering
	\includegraphics[width=.8\linewidth]{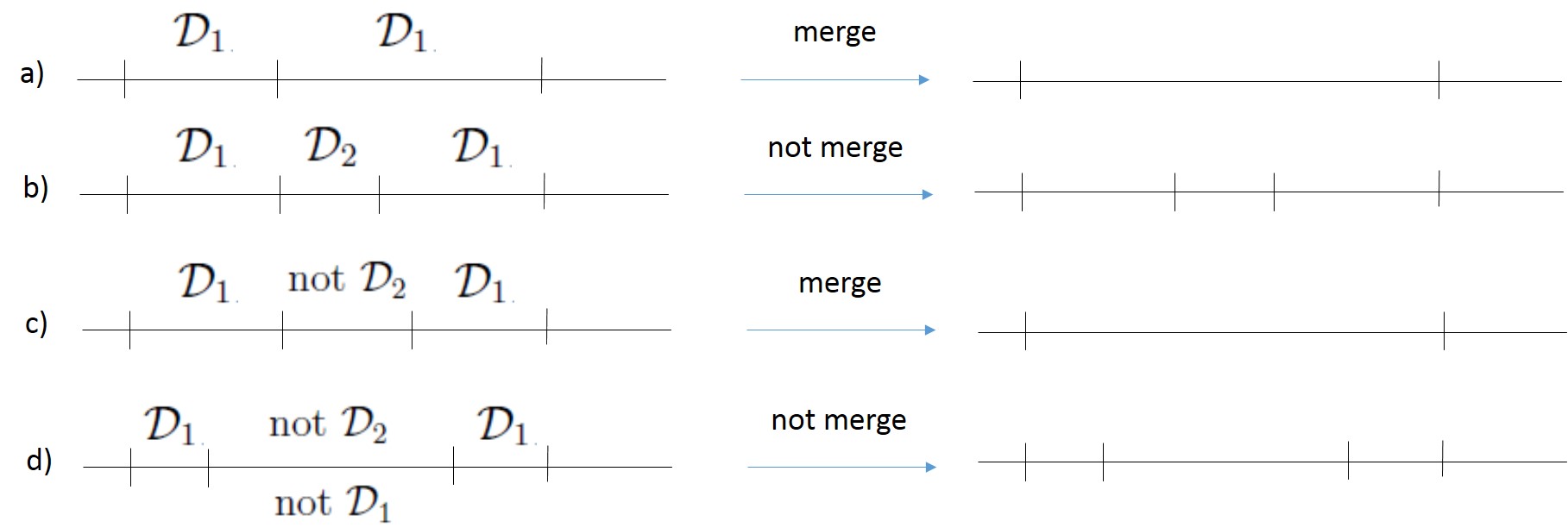}
	\caption{Illustration of the merge procedure for the intervals.}
	\label{fig::merge}
\end{figure*}

Then, the following sets are defined: $\mathcal{D}_1 = \{k: \estboldalpha_k^- 
\neq 0\}$ (called ``strong non zeros'') and $\mathcal{D}_2 = \{k: 
\estboldalpha_k^+ = 0\}$ (called ``strong zeros''). This step is illustrated in 
Figure \ref{fig::path}.
Intervals are merged using 
the following rules:
\begin{itemize}
	\item ``neighbor rule'': consecutive intervals of the same set are merged 
($\tau_k$ and $\tau_{k+1}$ are merged if both $k$ and $k+1$ belong to 
$\mathcal{D}_1$ or if they both belong to $\mathcal{D}_2$) (see a) and b) in 
Figure \ref{fig::merge});
	\item ``squeeze rule'': $\tau_k$, $\tau_{k+1}$ and $\tau_{k+2}$ are merged if 
both $k$ and 
$k+2$ belong to $\mathcal{D}_1$ while $k+1 \notin \mathcal{D}_2$ (or if both 
$k$ 
and $k+2$ belong to $\mathcal{D}_2$ while $k+1 \notin \mathcal{D}_1$) and $l_k 
+ l_{k+2} > l_{k+1}$ with $l_k = \max \tau_k - \min \tau_k$ (see c) and d) in 
Figure \ref{fig::merge}).
\end{itemize}

If the current value of $P$ does not yield any fusion between intervals, $P$ is 
updated by $P 
\leftarrow P + P_0$ in which $P_0$ is the initial value of $P$. 
The procedure is iterated until all the original intervals have been merged.

The result of the method is a collection of models 
$(\estboldalpha_k^*)_{k=1,\ldots,D}$, starting with $p$ intervals and finishing 
with one. The final selected model is the one that minimizes the CV error. In 
practice, this often results in a very small number of contiguous intervals 
which are of the same type (zero or non zero) and are easily interpretable (see 
Section~\ref{simulations}).

Let us remark that the intervals $(\tau_k)_{k=1,\ldots,D}$ are not used in 
the ridge step of Section~\ref{ilsir-ridge}, which can thus be performed once,
independently of the interval search. 
The whole procedure is described in Algorithm \ref{algo::algototal}.

\begin{algorithm}[ht]
	\caption{Overview of the complete procedure}
	\label{algo::algototal}
	\begin{algorithmic}[1]
		\State {\bf Ridge estimation} 
			\State Choose $\mu_2$ and $d$ according to Section \ref{param-selection}
			\State Solve the ridge penalized SIR to obtain $\hat{A}$ and $\hat{C}$, 
ridge estimates of the SIR (see details in Section \ref{ilsir-ridge})
		\State {\bf Sparse estimation} 
		\State Initialize the intervals $(\tau_k)_{k=1,\ldots,D}$ to $\tau_k = 
\{t_k\}$
		  \Repeat 
			\State Estimate and store $(\estboldalpha_k^*)_{k=1,\ldots,D}$ the 
solutions of the sparse 
problem that minimizes the GCV error
	\State Estimate $(\estboldalpha_k^+)_{k=1,\ldots,D}$ and 
$(\estboldalpha_k^-)_{k=1,\ldots,D}$  such that at most (resp. at least) a 
proportion $P$ of the 
coefficients are non zero coefficients (resp. are zero coefficients), for a 
given chosen $P$ (details in Section \ref{interv})
			\State Update  the intervals $(\tau_k)_{k=1,\ldots,D}$ according to the 
``neighbor'' and the ``squeeze'' rules (see Section \ref{interv})
			\Until{$\tau_1 \ne \left[t_1,t_p\right]$}
			\State Output : A collection of models 
$(\estboldalpha_k^*)_{k=1,\ldots,D}$
\State Select the model $(\estboldalpha_k^*)^*_{k=1,\ldots,D^*}$ that minimizes 
the CV error
\State Active intervals (for interpretation) are consecutive $\tau_k$ with non 
zero coefficients $\estalpha^*_k$
	\end{algorithmic}
\end{algorithm}

\section{Choice of parameters in the high dimensional setting}
\label{param-selection}

The method requires to tune four parameters : the number of slices $H$, the 
dimension of the EDR space $p$, the penalization parameter of the ridge 
regression $\mu_2$ and of the one of the sparse procedure $\mu_1$. Two of these 
parameters, $H$ and $\mu_1$, are chosen in a standard way (see 
Section~\ref{simulations} for further details).
This section presents a method to jointly choose $\mu_2$ and $d$,
for which no solution has been proposed that is suited to our high-dimensional framework.
Two issues are raised to tune these two parameters: i) they depend from each 
other and 
ii) the existing methods to tune them are only valid in a low-dimensional 
setting ($p < n$). 
We propose an iterative method which adapts existing approaches only valid for 
the low 
dimension framework and combine them to find an optimal joint choice for 
$\mu_2$ and $d$.

\subsection{A Cross-Validation (CV) criterion for $\mu_2$}
\label{choice_mu2}

Using the results of \cite{golub_etal_T1979}, \cite{li_yin_B2008} propose a 
Generalized Cross-Validation (GCV) 
criterion to select the regularization parameter $\mu_2$ and 
\cite{bernardmichel_etal_B2008} explain that this criterion can be applied to 
their modified estimator, using similar calculations. However, it requires the 
computation of $\estSigma^{-1/2}$, which does not exist in the high dimensional 
setting. 

We thus used a different strategy, based on $L$-fold cross-validation (CV), 
which is 
also used to select the best dimension of the EDR space, $d$ (see 
Section~\ref{choice_d}). More precisely, the data are split into $L$ folds, 
$\mathcal{L}_1$, \ldots, $\mathcal{L}_L$ and a CV error is 
computed for several values of $\mu_2$ in a given search grid and for a given 
(large enough $d_0$). The optimal $\mu_2$ is chosen as the one minimizing the 
CV error for $d_0$. 

The CV error is computed based on 
the original regression problem $\mathcal{E}_1(A,C)$. In the expression of 
$\mathcal{E}_1(A,C)$ and for the iteration number $l$ ($\in\{1,\ldots,L\}$), 
$A$ and $C_h$ are replaced by their estimates computed without the observations 
in fold number $l$. Then, an error is computed by replacing the values of 
$\hat{p}_h$, $\meanX_h$, $\meanX$ and $\estSigma$ by their empirical estimators 
for the observations in fold $l$. The precise expression is given in step 5 of 
Algorithm ~\ref{algo::choice-parameters} in Appendix~\ref{app::hyperparameters}.

\subsection{Choosing $d$ in a high dimensional setting}
\label{choice_d}

The results of CV ({\it i.e.}, the values of $\mathcal{E}_1(A,C)$ estimated by 
$L$-fold CV) are not directly usable for tuning $d$. The reason is similar to 
the one developed in \citet{biau_etal_IEEETIT2005,fromont_tuleau_COLT2006}: 
different $d$ correspond to different MLR (Multiple Linear Regression) problems 
which cannot be compared 
directly using a CV error. In such cases, an additional penalty depending on 
$d$ is necessary to perform a relevant selection and avoid overfitting due to 
large $d$.

Alternatively, a number of works have been dealing with the choice of $d$ in 
SIR. Many of them are asymptotic methods \citep{li_JASA1991,schott_JASA1994,
bura_cook_JASA2001,cook_yin_ANZJS2001,bura_yang_JMA2011,liquet_saracco_CS2012}
which are not directly applicable in the high dimensional framework. When $n < 
p$, \citet{zhu_etal_JASA2006,li_yin_B2008} estimate $d$ using the number of non 
zero eigenvalues of $\Gamma$, but their approach requires setting a 
hyper-parameter to which the choice of $d$ is sensitive. 
\cite{portier_delyon_JASA2014} describes an efficient approach that can be used 
when $n < p$ but it is based on bootstrap sampling and would thus be overly extensive 
in our situations where $d$ has to be tuned jointly with $\mu_2$ (see next 
section). 

Another point of view can be taken from \citet{li_JASA1991} who introduces a 
quantity, denoted by $R^2(d)$, which is the average of the squared canonical 
correlation between the space spanned by the columns of $\Sigma^{1/2} A$ and the 
columns of the space spanned by the columns of $\estSigma^{1/2} \widehat{A}$. As 
explained in \citet{ferre_JASA1998}, a relevant measure of quality for the 
choice of a dimension $d$ is $R(d) = d - \mathbb{E}\left[\mbox{Tr}\left(\Pi_d 
\estPi_d\right)\right]$, in which $\Pi_d$ is the $\Sigma$-orthogonal projector 
onto the subspace spanned by the columns of $A$ and $\estPi_d$ is the 
$\estSigma$-orthogonal projector onto the space spanned by the columns of 
$\hat{A}$. This quantity is equal to $\frac{1}{2} \mathbb{E} \left\| \Pi_d - 
\estPi_d \right\|^2_F$ (in which $\|.\|_F$ is the Frobenius norm; see the proof 
in Appendix~\ref{app:R2d}).

In practice, the quantity $\Pi_d$ is unknown and 
$\mathbb{E}\left[\mbox{Tr}\left(\Pi_d \estPi_d\right)\right]$ is thus 
frequently estimated by re-sampling techniques as bootstrap. Here, we choose a 
less computationally demanding approach by performing a CV estimation: 
$\mathbb{E}\left[\mbox{Tr}\left(\Pi_d \estPi_d\right)\right]$ is estimated 
during the same $L$-fold loop described in Section~\ref{choice_mu2}. An 
additional problem comes from the fact that, in the high dimensional setting, 
the $\estSigma$-orthogonal projector onto the space spanned by the columns of 
$\hat{A}$ is not well defined since the matrix $\estSigma$ is ill-conditioned.
This estimate is replaced by its regularized version using the 
$(\estSigma+\mu_2\mathbb{I}_p)$-orthogonal projector onto the space spanned by 
the columns of $\hat{A}$ and $\estPi_d$ is the 
$(\estSigma+\mu_2\mathbb{I}_p)$-orthogonal projector onto the space spanned by 
the columns of $\hat{A}$. Finally, for all $l=1,\ldots,L$, we computed the 
$(\estSigma^{\setminus l}+\mu_2\mathbb{I}_p)$-orthogonal projector onto the 
space spanned by the columns of $\hat{A}(l)$ in which $\estSigma^{\setminus l}$ 
and $\hat{A}(l)$ are computed without the observations in fold number $l$ and 
averaged the results to obtain an estimate of 
$\mathbb{E}\left[\mbox{Tr}\left(\Pi_d \estPi_d\right)\right]$. 

In practice, this estimate is often a strictly increasing function of $d$ and 
we chose the optimal dimension as the largest one before a gap in this 
increase (``elbow rule'').

\subsection{Joint tuning}

The estimation of $\mu_2$ and $d$ is jointly performed using a single CV pass 
in which both parameters are varied. Note that only the number of different 
values for $\mu_2$ strongly influences the computational time since SIR 
estimation is only performed once for all values of $d$, and selecting the 
first $d$ columns of $\widehat{A}$ for the last computation of the two 
criteria, the estimation of $\mathcal{E}(A,C)$ and that of $R(d)$. The 
overall method is described in Appendix~\ref{app::hyperparameters}.

\section{Experiments}\label{simulations}

This section evaluates different aspects of the methods on simulated and real 
datasets. The relevance of the selection procedure is evaluated on simulated 
and real datasets in Sections~\ref{simulated} and \ref{sec:sunflo}. 
Additionally, its efficiency in a regression framework is assessed on a real 
supervised regression problem in Section~\ref{tecator}. 

All experiments have been performed using the \RR{} package \pkg{SISIR}.  
Datasets and \RR{} scripts are provided at 
\url{https://github.com/tuxette/appliSISIR}.

\subsection{Simulated data}\label{simulated}

\subsubsection{Model description}
To illustrate our approach, we first consider two toy datasets, built as follow: 
$Y = \sum_{j=1}^d \log 
\left|\langle X, \bolda_j \rangle \right|$ with $X(t) = Z(t) + \epsilon$ in 
which $Z$ is a Gaussian process indexed on $[0,1]$ with mean $\mu(t) = -5 + 4t 
- 
4t^2$ and the Matern 3/2 covariance function 
\citep{rasmussen_williams_GPML2006}, and $\epsilon$ is a centered 
Gaussian variable independent of $Z$. The vectors $\bolda_j$ have a sinusoidal 
shape, but are nonzero only on specific intervals $I_j$: $\bolda_j = \sin 
\left(\frac{t(2+j)\pi}{2} - \frac{(j-1)\pi}{3}\right) \mathbb{I}_{I_j}(t)$. 

From this basis, we consider two models with increasing complexity:
\begin{itemize}
 \item {\bf (M1)}: $d=1$, $I_1 = [0.2,0.4]$
 \item {\bf (M2)}: $d=3$ and $I_1 = [0,0.1]$, $I_2=[0.5,0.65]$ and 
$I_3 = [0.65,0.78]$.
\end{itemize}

For both cases the datasets consist of $n=100$ observations of $Y$, digitized 
at $p = 200$ and $300$ evaluation points, respectively. The number of slices 
used to estimate the conditional mean $\mathbb{E}(X|Y)$ has been chosen equal 
to $H=10$: according to \citet{li_JASA1991,coudret_etal_JSFdS2014} among 
others, the performances of SIR estimates are not sensitive to the choice of 
$H$, as long as it is large enough (on a theoretical point of view, $H$ is 
required to be larger than $d+1$).

\begin{figure}[ht]
	\centering
	\includegraphics[width=\linewidth]{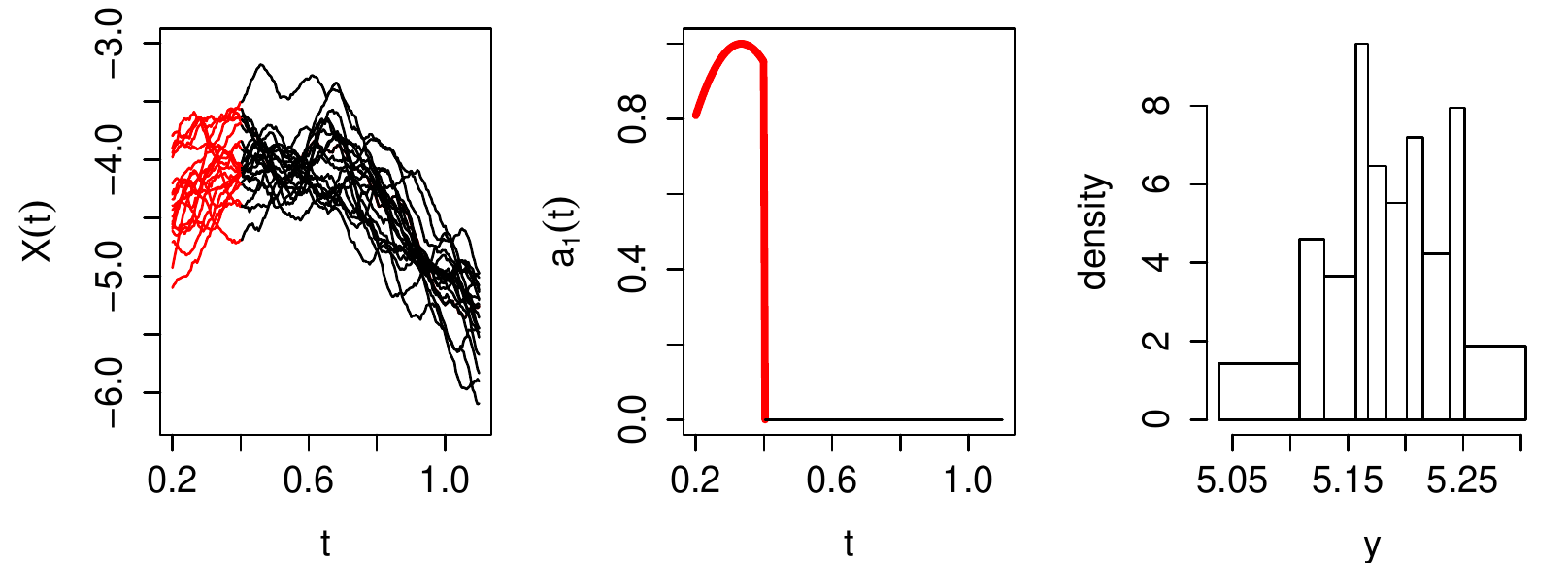} 
	\includegraphics[width=\linewidth]{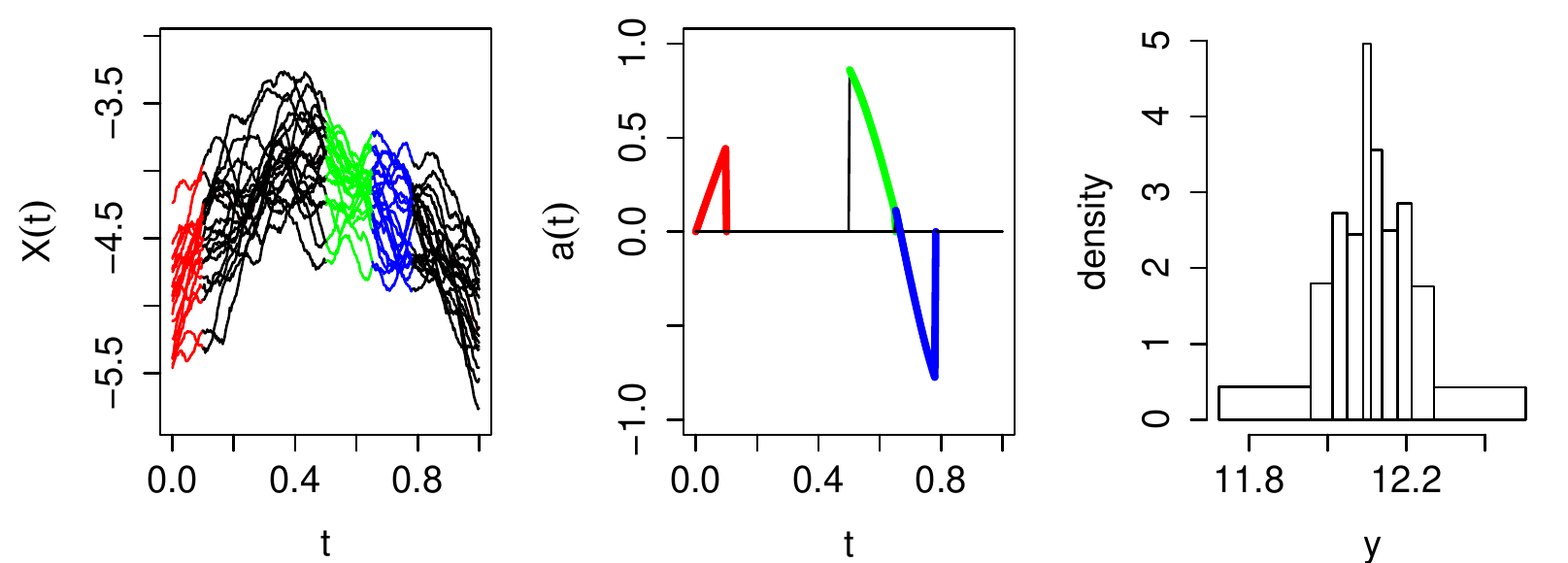}
		\caption{Summary of the two simulated datasets: top {\bf (M1)}, bottom 
{\bf (M2)}. The left charts display ten samples of $X$, the colors showing
the actual relevant intervals; the middle charts display the 
functions that span the EDR space with the relevant slices highlighted in 
color; the right charts display the distribution of $Y$.}
	\label{fig::simulated-data}
\end{figure}

The datasets are displayed in 
Figure~\ref{fig::simulated-data}, with {\it a priori} intervals provided to 
test the sparse penalty (see Section~\ref{simulated-sparse} for further 
details).

\subsubsection{Step 1: Ridge estimation and parameter selection}

\begin{figure}[ht]
	\centering
	\begin{tabular}{cc}
		\includegraphics[width=0.45\linewidth]{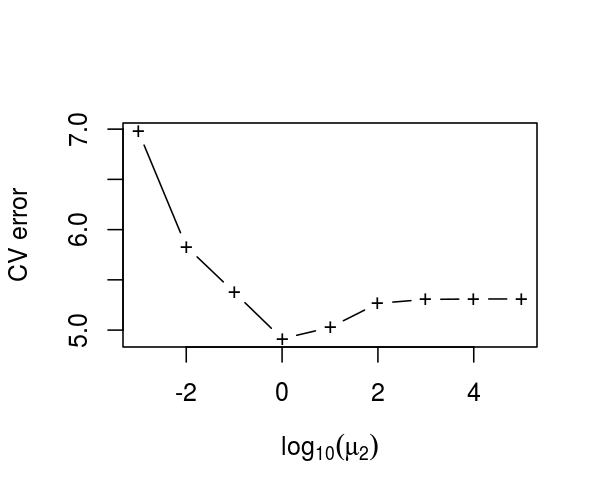} &
		\includegraphics[width=0.45\linewidth]{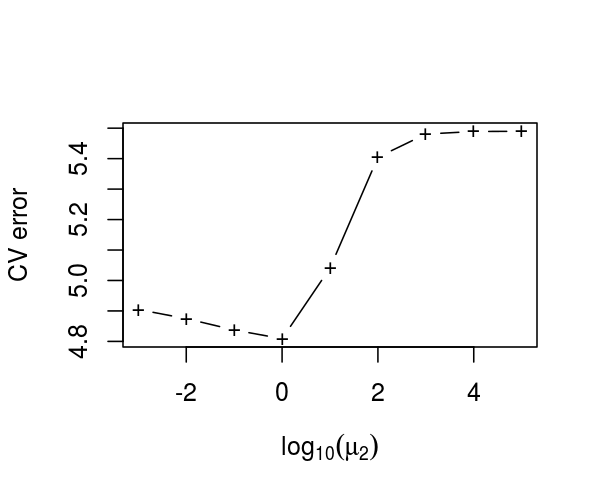}\\
		\includegraphics[width=0.45\linewidth]{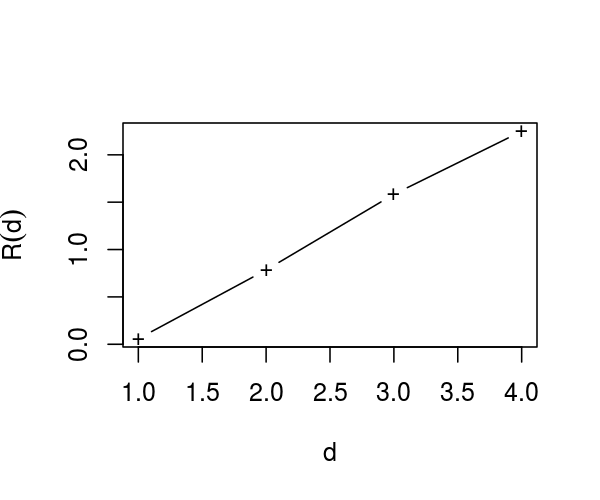} &
		\includegraphics[width=0.45\linewidth]{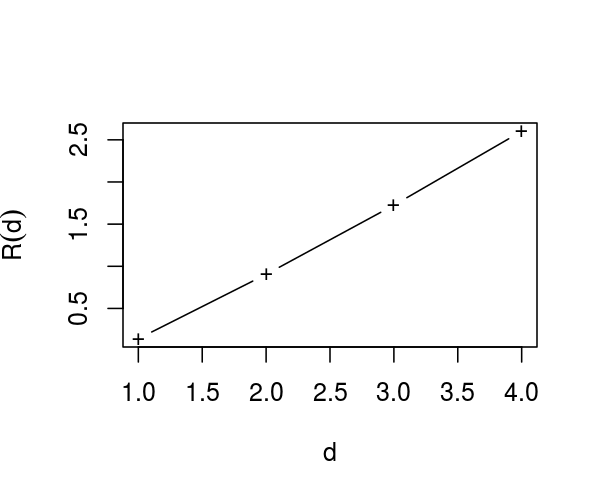}\\
		{\bf (M1)} & {\bf (M2)}
	\end{tabular}
	\caption{Top: CV error versus $\mu_2$ ($\log_{10}$ scale, for $d = 1$) and 
Bottom: estimation of $\mathbb{E}(R(d))$ versus $d$ (for $\mu_2$ = 1 in both 
cases), for models {\bf (M1)} (left) and {\bf (M2)} (right).}
	\label{fig::cverror}
\end{figure}

The method described in Section~\ref{ilsir-ridge} with parameter selection as 
in Section~\ref{param-selection} has been used to obtain the ridge estimates of 
$(\bolda_j)$ and to select the parameters $\mu_2$ (ridge regularization) and 
$d$ (dimension of the EDR space). Figure~\ref{fig::cverror} shows the evolution 
of the CV error and of the estimation of $\mathbb{E}(R(d))$ versus 
(respectively) $\mu_2$ and $d$ among a grid search both for $\mu_2 \in 
\{10^{-2}, 10^{-1}, ..., 10^5\}$ and $d \in \{1, 2, \ldots, 10\}$.
The chosen value for $\mu_2$ is 1 for both models and the chosen values for $d$, 
given by the ``elbow rule'' are $d=1$ for both 
models. The true values are, respectively, $d=1$ and $d=3$, which shows that 
the criterion tends to slightly underestimate the model dimension.

\subsubsection{Step 2: Sparse selection and definition of relevant intervals}
\label{simulated-sparse}

\begin{figure}[ht]
	\centering
	\begin{tabular}{cc}
		\includegraphics[width=0.4\linewidth]{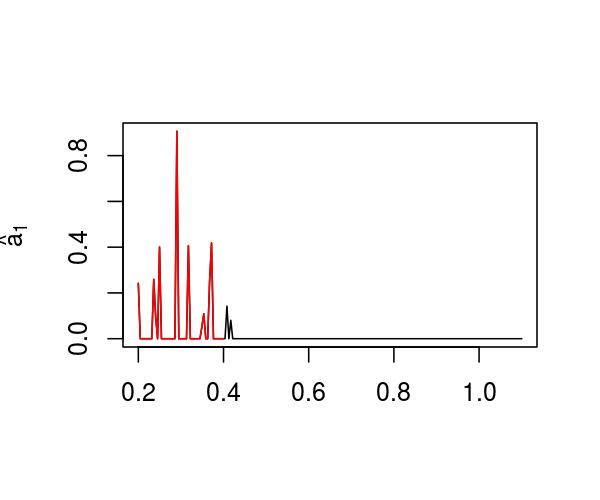} &
		\includegraphics[width=0.4\linewidth]{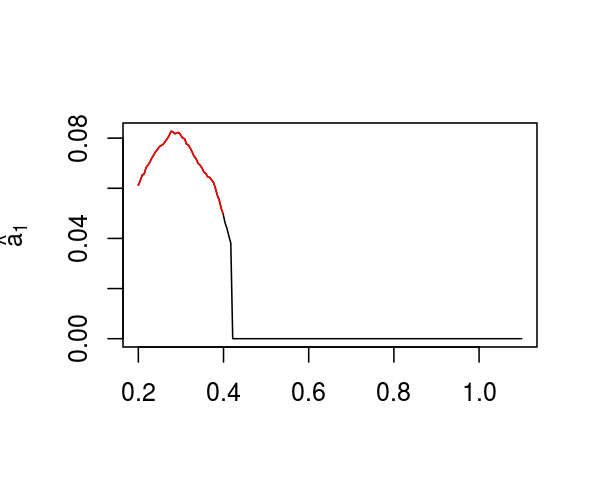} \\
		(a) 200 intervals &(b) 142 intervals\\
		\includegraphics[width=0.4\linewidth]{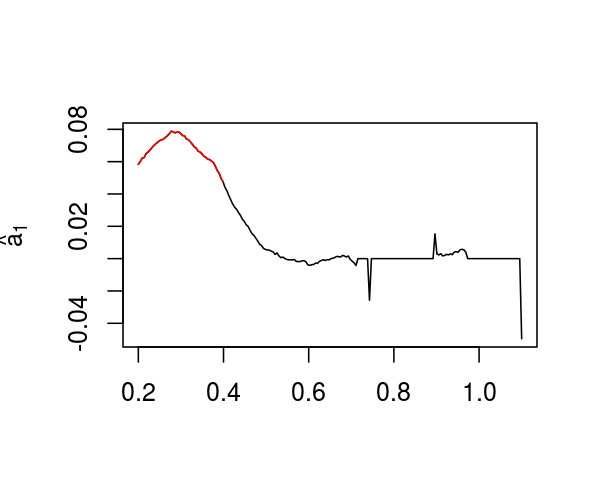}&
		\includegraphics[width=0.4\linewidth]{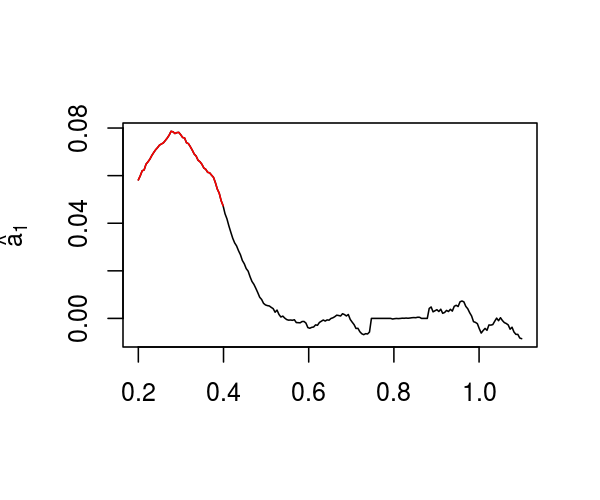}\\
		(c) 41 intervals & (d) 5 intervals\\
	\end{tabular}
	\caption{{\bf (M1)} Values of $\hat{\bolda}_1^s$ corresponding to four models 
	obtained using the iterative procedure with a different numbers of intervals.
(b) is the chosen model and (a) corresponds to a standard sparse estimation 
with no constraint on intervals.}
	\label{figmodel}
\end{figure}

The approach described in Section~\ref{interv} is then applied to both models.
The algorithm produces a large collection of models with a 
decreasing number of intervals: a selection of the estimates of $\bolda_1$ 
for {\bf (M1)},
corresponding to those models is shown in Figure~\ref{figmodel}.
The first chart (Figure \ref{figmodel},a) corresponds to the standard 
sparse penalty in which the constraint is put on isolated evaluation points. 
Even though most of selected points are found in the relevant interval, the estimated parameter 
$\hat{\bolda}_1^s$ has an uneven aspect which does not favor interpretation. 

By constrast, for a low number of intervals (less than 50, Figure 
\ref{figmodel}, c and d), the selected
relevant points (those corresponding to nonzero coefficients) have a much 
larger range than the original relevant interval (in red on the figure). 

The 
model selected by minimization of the cross-validation error (Figure \ref{figmodel}, b) was found 
relevant: this approach lead us to choose the model with 142 intervals, 
which actually correspond to two distinct and 
consecutive intervals (a first one, which contains only nonzero coefficients 
and a second one in which no point is selected by the sparse estimation). This 
final estimation is very close to the actual direction $\bolda_1$, both in 
terms of shape and support.

\begin{figure}[ht]
	\centering
	\begin{tabular}{cc}
		\includegraphics[width=0.55\linewidth]{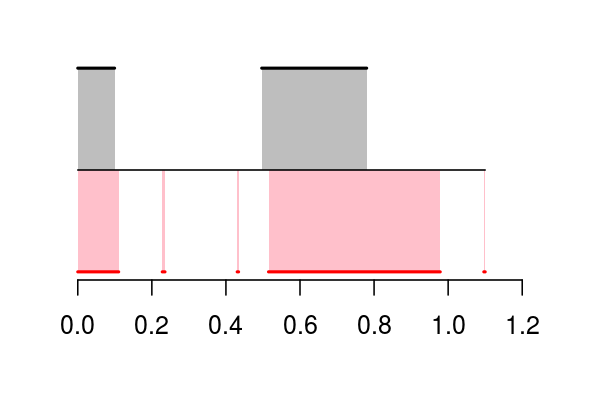} &
		\includegraphics[width=0.4\linewidth]{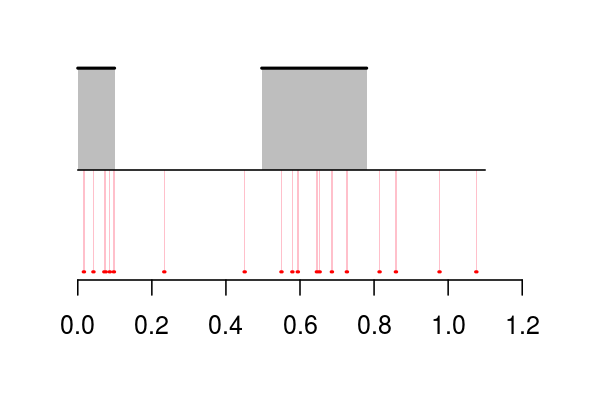} \\
		SISIR & standard sparse
	\end{tabular}
	\caption{{\bf (M2)} Left: comparison between the true intervals and the 
estimated ones. True intervals are represented in the upper side of the 
figure (in black) and by the gray background. Estimated intervals are 
represented by the red lines in the bottom of the figure and by the pink 
background. Right: same representation for the standard sparse approach 
(penalty is applied to $t_j$ and not to the intervals).}
	\label{fig::figM2}
\end{figure}

The same method is used for {\bf(M2)}. A comparison between the true relevant 
intervals and the estimated ones is provided in Figure \ref{fig::figM2} (left).
The support of each of the estimate $\hat{\bolda}_1$ is fairly appropriate: it 
slightly overestimates the length of the two real intervals and contains only 
three additional isolated points which are not relevant for the estimation. 
Compared to the standard sparse approach (right part of 
Figure~\ref{fig::figM2}), the approach is much more efficient to select the 
relevant intervals and provide more interpretable results by identifying 
properly important contiguous areas in the support of the predictors.

\begin{figure}[ht]
	\centering
	\begin{tabular}{cc}
		\includegraphics[width=0.45\linewidth]{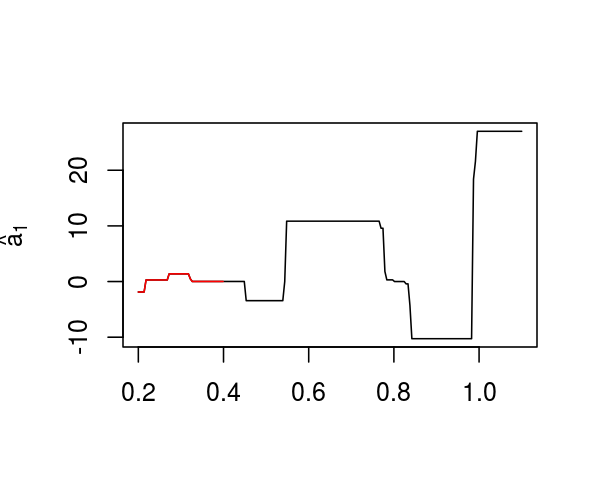} &
		\includegraphics[width=0.5\linewidth]{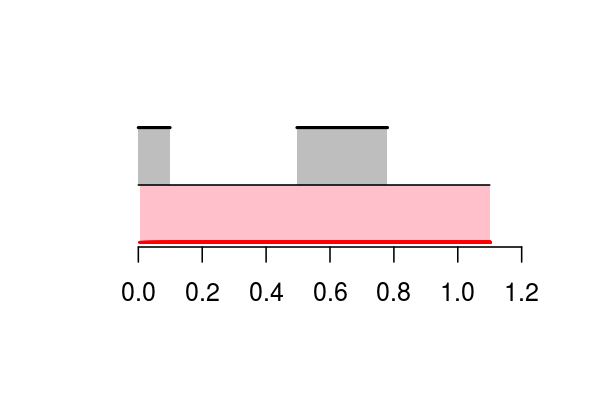}\\
		{\bf (M1)} & {\bf (M2)}
	\end{tabular}
	\caption{{\bf (M1)} Values of $\hat{\bolda}_1^s$ obtained with fused Lasso. 
The target relevant interval is highlighted in red. {\bf (M2)} Comparison 
between the true intervals and the estimated ones. True intervals are 
represented in the upper side of the figure (in black) and by the gray 
background. Fused Lasso estimated intervals are represented by the red lines in 
the bottom of the figure and by the pink background.}
	\label{fig::simulated-fused}
\end{figure}

As a basis for comparison, fused Lasso \citep{tibshirani_etal_JRSSB2005}, as 
implemented in the \RR{} package \pkg{genlasso}, was used with both {\bf (M1)} 
and {\bf (M2)} datasets. For comparison with our method, we applied fused Lasso 
on the output of the ridge SIR so as to find $\bolda_1^s \in 
\R{p}$ that minimizes:
\begin{eqnarray*}
	&&\sum_{i=1}^n \left[\mathcal{P}_{\hat{\bolda}_1}(X|y_i) - (\bolda_1^s)^\top 
\x_i \right]^2 +\\
	&& \qquad \lambda_1 \|\bolda_1^s\|_{\ell_1} + \lambda_2 \sum_{j=1}^{p-1} 
|a_{1j}^s - a^s_{1,j+1}|,
\end{eqnarray*}
for $\bolda^s_1 = (a^s_{1,1}, \ldots, a^s_{1,p})$. The tuning parameters $\lambda_1$ and $\lambda_2$ 
were selected by 10-fold CV over a 2-dimensional grid search. The idea behind 
fused Lasso is to have a large number of identical consecutive entries in 
$\bolda_1^s$. In our framework, the hope is to automatically design relevant 
intervals using this property. Results are displayed in 
Figure~\ref{fig::simulated-fused} for both simulated datasets.
Contrary to simple Lasso, fused Lasso produces a piecewise constant estimate. 
However, both for {\bf (M1)} and {\bf (M2)}, the method fails to provide a 
sparse solution: almost the whole definition domain of the predictor is 
returned as relevant.

\subsection{Tecator dataset}\label{tecator}

\begin{figure}[ht]
	\centering
	\includegraphics[width=0.9 \linewidth]{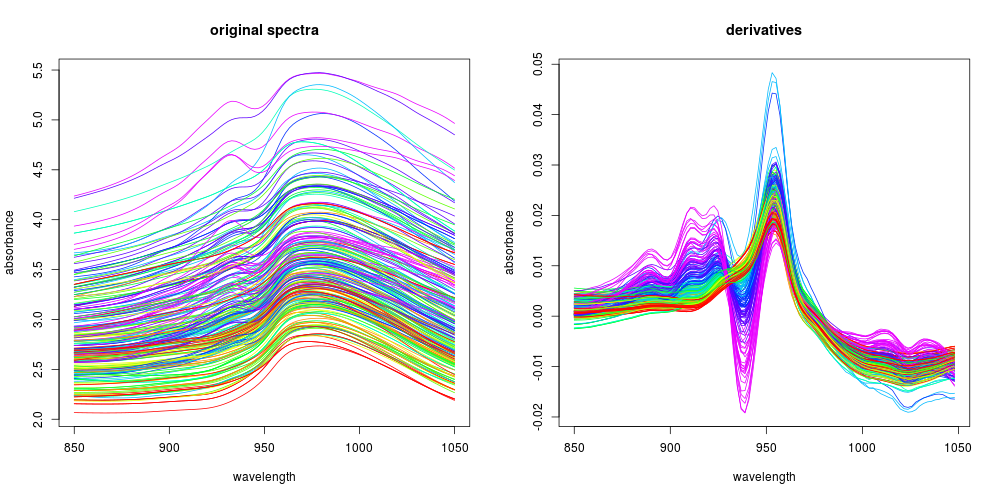}
	\caption{\textbf{Tecator}. 215 near infrared spectra from the ``Tecator'' 
dataset (left) and corresponding first order derivatives (right). The color 
level indicates the percentage of fat content.}
	\label{fig::tecator}
\end{figure}

Additionally, we tested the approach with the well-known Tecator dataset
\citep{borggaard_thodberg_AC1992}, which consists of spectrometric data from 
the food industry. This dataset is a standard benchmark for functional data 
analysis. It contains 215 observations of near infrared absorbency spectra of a 
meat sample recorded on a Tecator Infratec Food and Feed Analyzer. Each 
spectrum was sampled at 100 wavelengths uniformly spaced in the range 
850--1050~nm. The composition of each meat sample was determined by analytic 
chemistry, among which we focus on the percentage of fat content. The data is 
displayed in Figure~\ref{fig::tecator}: the left chart displays 
the original spectra whereas the right chart displays the first order 
derivatives (obtained by simple finite differences). The fat content is 
represented in both graphics by the color level and, as is already well known 
with this dataset, the derivative is a good predictor of this quantity: these 
derivatives were thus used as predictors ($X$) to explain the fat content ($Y$).

We first applied the method on the entire dataset to check the relevance 
of the estimated EDR space and corresponding intervals in the domain 
850--1050~nm. Using the ridge estimation and the method described in 
Section~\ref{param-selection}, we set $\mu_2 = 10^{-4}$ and $d=1$. 

The relevance of the approach was then assessed in a regression setting. 
Following the simulation setting described in \cite{hernandez_etal_SS2015}, we 
split the data into a training and test sets with 150 observations for the 
training. This separation of the data was performed 100 times randomly. For each
training data set, the EDR space was estimated and the projection of the 
predictors on this space obtained. A a Support Vector Machine \citep[SVM, 
$\epsilon$-regression method, package \pkg{e1071}][]{meyer_etal_e10712015} was 
used to fit the link function $F$ of Equation~(\ref{eq::sir-model}) from both 
the projection on the EDR space obtained by a simple ridge SIR and the 
projection on the EDR space obtained by SISIR. The mean square error was then 
computed on the test set. We found an averaged value equal to 5.54 for the 
estimation of the EDR space obtained by SISIR and equal to 11.11 when the 
estimation of the EDR space is directly obtained by ridge SIR only. The 
performance of SISIR in this simulation is thus half the value reported for the 
Nadaraya-Watson kernel estimate in \cite{hernandez_etal_SS2015}. 

\begin{figure}[ht]
	\centering
	\includegraphics[width=0.45\linewidth]{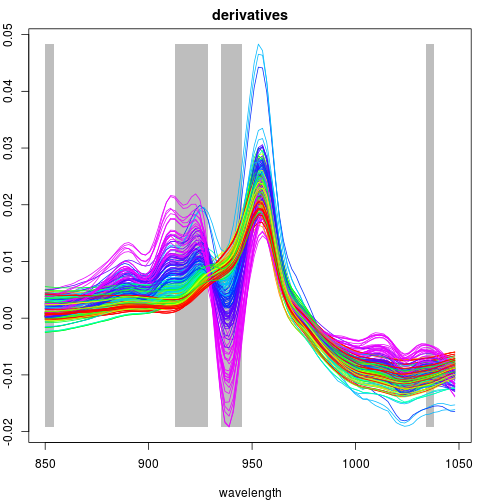} 
	\includegraphics[width=0.4\linewidth]{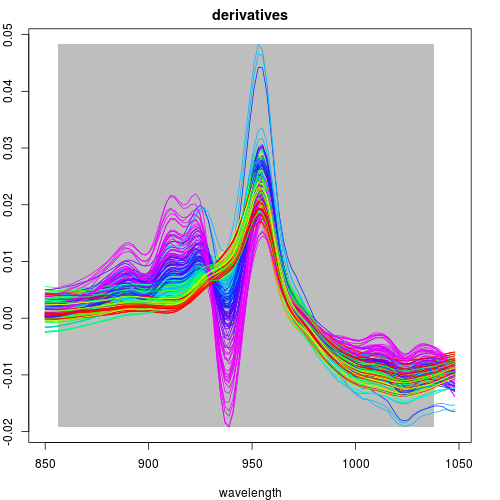}
	\caption{\textbf{Tecator}. Left: original predictors (first order 
derivatives) 
with a gray background superimposed to highlight the active 
intervals found by our procedure. Right: same figure for the standard sparse 
approach (no constraint on intervals).}
	\label{fig::tecator-esta}
\end{figure}

Even if some methods achieve better performance on this data set 
(\cite{hernandez_etal_SS2015} reported an average MSE of 2.41 for their non 
parametric approach), our method has the advantage of being easily interpretable 
because it extracts a few components which are themselves composed of a small 
number of relevant intervals: Figure~\ref{fig::tecator-esta} shows the 
intervals selected in the simulation with the smallest MSE, compared to the 
values selected by the standard Lasso. Our method is able to identify two 
intervals in the middle of the wavelength definition domain that are actually 
relevant to 
predict the fat content (according to the ordering of the colors in this area). 
On the contrary, standard sparse SIR selects almost the entire interval. 

\subsection{Sunflower yield}\label{sec:sunflo}

Finally, we applied our strategy to a challenging agronomic problem,
the inference of interpretable climate-yield relationships on complex crop 
models.

We consider a process-based crop model called SUNFLO, which was developed to 
simulate the annual grain yield 
(in tons per hectare) of sunflower cultivars, as a function of time, environment 
(soil and climate), 
management practice and genetic diversity \citep{casadebaig2011sunflo}. 
SUNFLO requires functional inputs in the form of climatic series. These series 
consist of daily measures of five
variables over a year: minimal temperature, maximal temperature, global 
incident radiation, precipitations and evapotranspiration. 

The daily crop dry biomass growth rate is calculated as an ordinary differential
equation function of incident photosynthetically active radiation, light 
interception efficiency and radiation
use efficiency.  Broad scale processes of this framework, the dynamics of leaf 
area, photosynthesis and biomass
allocation to grains were split into  finer processes (e.g leaf expansion and 
senescence, response functions to
environmental stresses).  Globally, the SUNFLO crop model has about 50 equations 
and 64 parameters (43 plant-related
traits and 21 environment-related). Thus, due to the complexity of 
plant-climate interactions and the strongly irregular nature of climatic data, 
understanding the relation between yield and climate is a particularly 
challenging task.

\begin{figure}[ht]
	\centering
	\includegraphics[width=\linewidth]{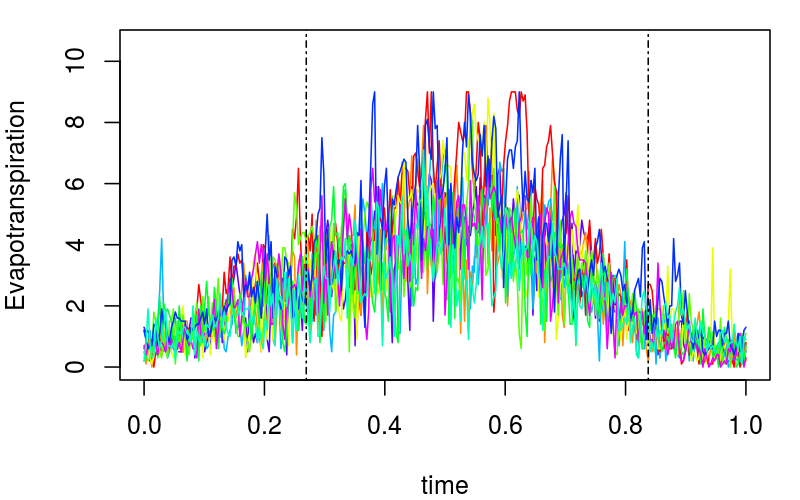} 
	\caption{\textbf{Sunflo}. Ten series of evaportranspiration daily recordings. 
The color level 
indicates the corresponding yield and the dashed lines bound the actual 
simulation definition domain.}
	\label{fig:etp}
\end{figure}

The dataset used in the experiment consisted of 111 yield values computed using
SUNFLO for different climatic series (recorded between
1975 and 2012 at five French locations). We focused 
solely on evapotranspiration  as a functional predictor
because it is essentially a combination of the other four variables 
\citep{allen1998crop}. 
The cultural year ({\it i.e.}, the period on which the simulation is performed) 
is from weeks 16 to 41 (April to October).
We voluntarily kept unnecessary data (11 weeks before simulation and 8 weeks 
after) for testing purpose (because these periods are known to be irrelevant 
for the prediction). The resulting curves contained 309 measurement points.
Ten series of this dataset are shown in Figure \ref{fig:etp}, with colors 
corresponding to the yield that we intend to explain: no clear relationship can 
be identified between the the value of the curves at any measurement point and 
the yield value.

Using the ridge estimation and the method described in 
Section~\ref{param-selection}, we set $\mu_2 = 10^{3}$ and $d=2$. 
Then, we followed the approach described in Section~\ref{interv} to design the 
relevant intervals.

\begin{figure}[ht]
	\centering
	\includegraphics[width=.9\linewidth]{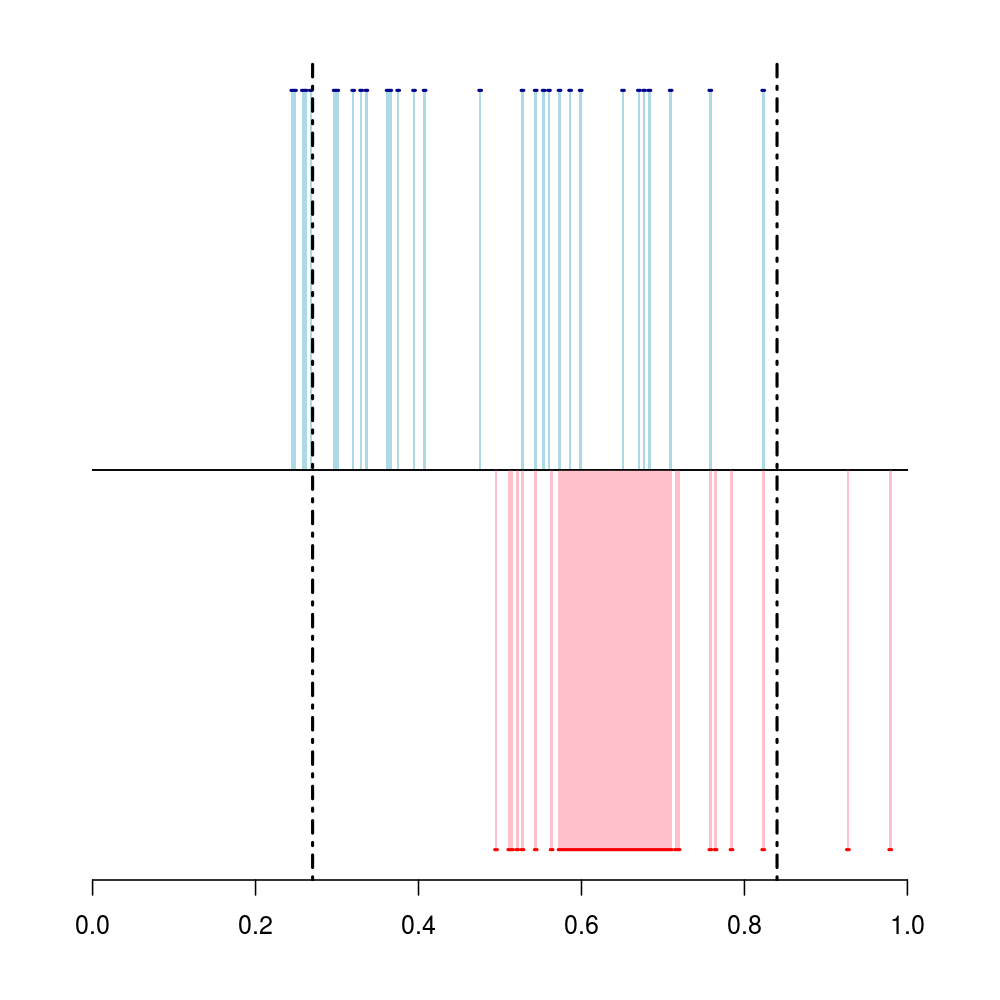} 
	\caption{\textbf{Sunflo}. Top: standard sparse SIR (blue). Bottom: SISIR 
(pink). The colored areas depict the active intervals. The dashed lines bound 
the actual simulation definition domain.}
	\label{fig:sunfloresults}
\end{figure}

Figure~\ref{fig:sunfloresults} shows the selected intervals obtained after 
running our algorithm, as well as
the points selected using a standard sparse approach. The standard sparse SIR 
(top of the figure) captures well the simulation interval 
(with only two points selected outside of it), but fails to identify the 
important periods within it.
In contrast, SISIR (bottom) focuses on the second half of the simulation 
interval, and in particular its third quarter.
This matches well expert knowledge, that reports little influence of the climate 
conditions at early stage of the plant growth and almost none once the grains 
are ripe \citep{casadebaig2011sunflo}.

\section*{Acknowledgments}

The authors thank the two anonymous referees for relevant
remarks and constructive comments on a previous version of the paper.

\appendix
\section{Equivalent expressions for $R^2(d)$}
\label{app:R2d}

In this section, we show that $R^2(d) = \frac{1}{2}\mathbb{E} \left\| \Pi_d - 
\widehat{\Pi}_d \right\|^2_F$. We have

\begin{eqnarray*}
	\frac{1}{2} \left\| \Pi_d - \widehat{\Pi}_d \right\|^2_F &=& \frac{1}{2} 
\mbox{Tr}\left[ \left( \Pi_d - \widehat{\Pi}_d \right) \left( \Pi_d - 
\widehat{\Pi}_d \right)^\top \right]\\
	&=& \frac{1}{2} \mbox{Tr}\left[ \left( \Pi_d \Pi_d \right)\right] - 
\mbox{Tr}\left[ \left( \Pi_d \widehat{\Pi}_d \right)\right] + \\
	&& \qquad \frac{1}{2} \mbox{Tr}\left[ \left( \widehat{\Pi}_d \widehat{\Pi}_d 
\right)\right].
\end{eqnarray*}

The norm of a $M$-orthogonal projector onto a space of dimension $d$ is 
equal to $d$, we thus have that
\[
	\frac{1}{2} \left\| \Pi_d - \widehat{\Pi}_d \right\|^2_F = d - 
\mbox{Tr}\left[ 
\left( \Pi_d \widehat{\Pi}_d \right)\right],
\]
which concludes the proof.

\section{Joint choice of the parameters $\mu_2$ and $d$}
\label{app::hyperparameters}

Notations: \begin{itemize}
	\item $\mathcal{L}_l$ are observations in fold number $l$ and 
$\overline{\mathcal{L}_l}$ are the remaining observations;
	\item $\hat{A}(\mathcal{L}, \mu_2, d)$ and $\hat{C}(\mathcal{L}, \mu_2, d)$ 
are minimizers of the ridge regression problem restricted to observations $i 
\in \mathcal{L}$. Note that for $d_1 < d_2$, $\hat{A}(\tau, \mu_2, d_1)$ are 
the first $d_1$ columns of $\hat{A}(\mathcal{L}, \mu_2, d_2)$ (and similarly 
for $\hat{C}(\mathcal{L}, \mu_2, d)$);
	\item $\hat{p}_h^\mathcal{L}$, $\meanX_h^\mathcal{L}$, $\meanX^\mathcal{L}$ 
and $\estSigma^\mathcal{L}$ are, respectively, slices frequencies, conditional 
mean of $X$ given the slices, mean of $X$ given the slices and covariance of 
$X$ 
for observations $i \in \mathcal{L}$;
	\item $\estPi_{d,\mu_2}^{\mathcal{L}}$ is the 
$(\estSigma^{\mathcal{L}}+\mu_2\mathbb{I}_p)$-orthogonal projector onto the 
space spanned by the first $d$ columns of $\hat{A}(\mathcal{L},\mu_2,d_0)$ and 
$\estPi_{d,\mu_2}$ is $\estPi_{d,\mu_2}^{\mathcal{L}}$ for $\mathcal{L} = 
\{1,\,\ldots,\,n\}$.
\end{itemize}

\begin{algorithm}[H]
	\caption{}
	\label{algo::choice-parameters}
	\begin{algorithmic}[1]
		\State Set $\mathcal{G}_{\mu_2}$ (finite search grid for $\mu_2$) and $d_0 
\in \mathbb{N}^*$ large enough
		\For{$\mu_2 \in \mathcal{G}_{\mu_2}$}
			\For{$l=1,\ldots,L$}
				\State Estimate $\hat{A}(\overline{\mathcal{L}_l},\mu_2,d_0)$ and 
$\hat{C}(\overline{\mathcal{L}_l},\mu_2,d_0)$
				\State With the observations $i\in\mathcal{L}_l$ and for $d\in \{1, 
\ldots, d_0\}$, compute 
				\begin{eqnarray*}
					\mbox{CVerr}^l_{d,\mu_2} &=& \sum_{h=1}^H \hat{p}_h^{\mathcal{L}_l} 
\left\|\left(\meanX_h^{\mathcal{L}_l} - \meanX^{\mathcal{L}_l} \right) - 
\right.\\
					&&\qquad \left. \estSigma^{\mathcal{L}_l} 
\hat{A}(\overline{\mathcal{L}_l},\mu_2,d)\hat{C}_h(\overline{\mathcal{L}_l},
\mu_2,d) \right\|^2_{(\estSigma^{\mathcal{L}_l}+\epsilon \mathbb{I})^{-1}}
				\end{eqnarray*}
				in which $\epsilon$ is a small positive number that makes 
$(\estSigma^{\mathcal{L}_l}+\epsilon \mathbb{I})$ invertible.
				\State For $d \in \{1, \ldots, d_0\}$, compute 
$\estPi_{d,\mu_2}^{\overline{\mathcal{L}_l}}$
				\State For $d \in \{1, \ldots, d_0\}$, compute 
				\[
					\hat{R}_{\mu_2}(d) = d - \frac{1}{L} \sum_{l=1}^L 
\mbox{Tr}\left(\estPi_{d,\mu_2}^{\overline{\mathcal{L}_l}} \estPi_{d,\mu_2} 
\right)
				\]
			\EndFor
			\State Compute 
			\[
				\mbox{CVerr}_{\mu_2,d} = \frac{1}{L} \sum_{l=1}^L 
\mbox{CVerr}^l_{\mu_2,d}
			\]
		\EndFor
		\State State $d^* \leftarrow d_0$.
		\Repeat
			\State Choose $\mu_2^* = \arg\min_{\mu_2\in\mathcal{G}_{\mu_2}} 
\mbox{CVerr}_{\mu_2,d^*}$
			\State Update $d^*$ with an ``elbow rule'' in $\hat{R}_{\mu_2^*}(d)$
		\Until{Stabilization of $d^*$}
		\State Output: $\mu_2^*$ and $d^*$
	\end{algorithmic}
\end{algorithm}

\bibliographystyle{spbasic}      
\bibliography{bibliototal}   

\end{document}